\documentclass[12pt]{article}

\usepackage{
	enumitem,
	times,
	amsmath,
	amssymb,
	amsthm,
	mathtools,
	fancyhdr,
	color,
	thmtools,
	etoolbox,
	xspace,
	hyperref,
	lscape,
	thm-restate,
	tikz-cd
}

\usepackage[margin=2cm]{geometry}

\declaretheorem[style=definition,numberwithin=section]{definition}
\declaretheorem[style=definition,qed=$\oslash$,sibling=definition]{example}

\declaretheorem[style=slanted,sibling=definition]{theorem}
\declaretheorem[style=slanted,sibling=definition]{conjecture}
\declaretheorem[style=slanted,sibling=definition]{lemma}

\declaretheorem[style=slanted,sibling=definition]{corollary}

\declaretheorem[style=definition,sibling=example]{remark}

\AtEndEnvironment{proof}{\setcounter{claim}{0}}

\usepackage[square, comma, numbers, sort&compress ]{natbib}
\newcommand{\C}{\ensuremath{\mathbb{C} }}
\newcommand{\R}{\ensuremath{\mathbb{R} }}
\newcommand{\N}{\ensuremath{\mathbb{N} }}
\newcommand{\Z}{\ensuremath{\mathbb{Z} }}

\DeclareMathOperator{\Ad}{Ad}

\newcommand\scalemath[2]{\scalebox{#1}{\mbox{\ensuremath{\displaystyle #2}}}}
\numberwithin{equation}{section}

\title{An addendum on the Mathieu Conjecture for $SU(N)$, $Sp(N)$ and $G_2$}
\author{Kevin Zwart\footnote{Email address: \href{mailto:kevin.zwart@ru.nl}{kevin.zwart@ru.nl}}}
\date{IMAPP, Radboud University Nijmegen, The Netherlands \\[2ex]}

\begin{document}
	
	\maketitle
	
	\begin{abstract}\noindent
		In this paper, we sharpen results obtained by the author in 2023. The new results reduce the Mathieu Conjecture on $SU(N)$ (formulated for all compact connected Lie groups by O. Mathieu in 1997) to a conjecture involving only functions on $\R^n\times (S^1)^m$ with $n,m$ non-negative integers instead of involving functions on $\R^n\times (S^1\setminus\{1\})^m$. The proofs rely on a more recent work of the author (2024) and a specific $KAK$ decomposition. Finally, with these results we can also improve the results on the groups $Sp(N)$ and $G_2$ in the latter paper, since they relied on the construction introduced in the 2023 paper.
	\end{abstract}
	\section{Introduction}
	
	In a recent paper \cite{Zwart}, we showed how the Mathieu conjecture for $SU(N)$ can be reduced to an abelian conjecture by using a generalized Euler Angles decomposition of $SU(N)$. However, in a more recent paper \cite{Zwart2}, we proved a theorem that allows a generalized Euler Angles decomposition of simply connected compact Lie groups. Since this proof is based on Lie theory instead of a case-by-case proof, the first part of the present paper applies mentioned theorem in \cite{Zwart2} to the group $SU(N)$, and reflect on the different results it gives with respect to \cite{Zwart2}.
	
	In addition, a more direct way of evaluating the integrals in Lemma 2.7 of \cite{Zwart} was found, which is given in Lemma \ref{lemma:rewrite_integral_for_SU(N)} below. The argument for this was actually hidden in the proof of Theorem 2.11 in \cite{Zwart}, and it makes Theorem 2.11 in \cite{Zwart} a corollary of Lemma \ref{lemma:rewrite_integral_for_SU(N)} below, see Theorem \ref{thm:Mathieu_proven_assuming_XZ-conjecture_SU(N)}. As a result, the resulting conjecture we ended with in \cite{Zwart} can be weakened. That is to say, in \cite{Zwart}, the functions that were involved in the conjecture were possible $N$-th roots of complex variables, while with the new approach, no $N$-th roots are necessary and one can just focus on Laurent polynomials. See for comparison Definition 2.9 and Conjecture 2.10 in \cite{Zwart} to Definition \ref{def:admissible_function} and Conjecture \ref{con:xz-conjecture} in this paper.
	
	In Section \ref{sec:Applying_to_Sp(N)_and_G_2}, we apply these results to $Sp(N)$ and $G_2$ as well, since they both required a decomposition of $SU(N)$ as well. As expected, we get a weaker conjecture for these groups as well, see Conjecture \ref{con:xz-conjecture_Sp(N)} and \ref{con:xz-conjecture_G_2}.
	
	\textit{Acknowledgments:} The author would like to thank Michael M\"uger for the continued help with the project. He also wishes to thank Erik Koelink for the insightful questions and the valued support.
	
	\section{Concerning the group $SU(N)$ and the Mathieu conjecture}\label{sec:applying_big_thm_to_SU(N)}
	In \cite{Zwart2}, we proved the following Theorem
	
	\begin{theorem}\emph{[The Euler angles Theorem]}\label{thm:KAK_decomp}
		Let $G$ be a simply connected compact Lie group, and let $\mathfrak{g}$ be its Lie algebra. Let $\theta:\mathfrak{g}\rightarrow\mathfrak{g}$ be an involutive automorphism. Let $\mathfrak{k}, \mathfrak{p}$ be the $+1$ and $-1$ eigenspace of $\theta$, respectively, in such a way that $\mathfrak{g}=\mathfrak{k}\oplus\mathfrak{p}$. Fix a maximal abelian subalgebra $\mathfrak{a}\subset \mathfrak{p}$, and let $\mathfrak{h}$ be any maximal abelian subalgebra containing $\mathfrak{a}$. Let $\Delta$ be the roots of  the complexification $\mathfrak{g}_\C$ with respect to $\mathfrak{h}_\C$, choose a set of positive roots $\Delta^+$, and define the set $\Delta_\mathfrak{p}:=\{\alpha\in \Delta|\,\,\alpha|_{\mathfrak{a}}\neq 0\}$. Furthermore, let $\Delta_\mathfrak{p}^+=\Delta^+\cap\Delta_\mathfrak{p}$ and define $\mathfrak{a}_+$ to be the closed positive Weyl chamber in $\mathfrak{a}$. Let $\mathcal{A}$ be the closure of the connected component of the set
		\begin{align*}
			\mathfrak{a}_{+}-\{H\in\mathfrak{a}\,|\,\alpha(H)\in \pi i\Z \,\text{ for some } \alpha\in \Delta_\mathfrak{p} \}
		\end{align*} in such a way that $0\in \mathcal{A}$. Finally let $K\subseteq G$ be the connected analytic Lie subgroup with Lie algebra $\mathfrak{k}$, and $M=Z_K(\mathfrak{a})$. Then the mapping 
		\begin{align}\label{eq:KAK_decomposition}
			(K/M)\times \exp(\mathcal{A})\times K \rightarrow G,\qquad\qquad (kM,\exp(H),l)\mapsto \exp(\Ad_G(k)H)l
		\end{align}
		is surjective, and a diffeomorphism up to a measure zero set if we replace $\mathcal{A}$ by $\mathrm{int}(\mathcal{A})$.
		
		In addition, define $J:A\rightarrow\C$ by
		\begin{align}\label{eq:Jacobian_abstract_Euler_decomp}
			J(\exp(H)) := \prod_{\alpha\in \Delta_\mathfrak{p}^+} \sin(\alpha(iH)).
		\end{align}
		Then the Haar measure decomposes in the parameterization given in (\ref{eq:KAK_decomposition}) as
		\begin{align*}
			\int_G h(g) dg = C \int_{K/M}\int_{\mathcal{A}}\int_K h(k_1\exp(H)k_2) |J(\exp(H))| dk_2 dH dg_{K/M}
		\end{align*}
		for any measurable function $h:G\rightarrow \C$, where $C>0$ is a constant (independent of $h$), $k_2\in K$ with corresponding Haar measure $dk_2$, $k_1\in k_1M$ an arbitrary representative of $k_1M\in K/M$ with corresponding unique $K$-invariant measure $dg_{K/M}$ on $K/M$, and $dH$ the measure on $\mathfrak{a}$.
	\end{theorem}
	
	In the first part of this paper we apply Theorem \ref{thm:KAK_decomp} to $SU(N)$ with $N\geq 2$. To do so, we define a spanning set of $\mathfrak{su}(N)$, as in \cite{Zwart,Zwart2}, given by 
	\begin{align*}
		[\lambda_{j^2-1+k}]_{\mu,\nu}&:=i(\delta_{\lceil\frac{k}{2}\rceil,\mu}\delta_{j+1,\nu}+\delta_{j+1,\mu}\delta_{\lceil\frac{k}{2}\rceil,\nu})\quad\qquad\text{if $k$ is odd},\\	[\lambda_{j^2-1+k}]_{\mu,\nu}&:=\delta_{\frac{k}{2},\mu}\delta_{j+1,\nu}-\delta_{j+1,\mu}\delta_{\frac{k}{2},\nu}\qquad\qquad\quad\text{if $k$ is even},\\ 
		[\lambda_{(j+1)^2-1}]_{\mu,\nu}&:= i(\delta_{j,\mu}\delta_{j,\nu}-\delta_{j+1,\mu}\delta_{j+1,\nu}),
	\end{align*}
	where $j=1,\ldots, N-1$ and $k=1,2,\ldots,2j$. In a similar way as in \cite{Zwart}, we get the following:
	
	\begin{lemma}\label{lemma:Euler_Angles_SU(N)_revisited}Let $N\geq 2$. Define inductively the map $F_{SU(N)}:[0,\pi)^{\frac{N(N-1)}{2}}\times \left(0,\frac{\pi}{2}\right)^{\frac{N(N-1)}{2}}\times[0,2\pi)^{N-1}\rightarrow SU(N)$ by $F_{SU(1)}\equiv 1$ and by
		\begin{equation}\label{eq:SU(N)_parametrization_new_version}
			\begin{split}
				&F_{SU(N)}(\phi_1,\ldots\phi_{\frac{N(N-1)}{2}},\psi_1,\ldots,\psi_{\frac{N(N-1)}{2}},\omega_1,\ldots,\omega_{N-1}):=\\
				&\scalemath{0.9}{\left(\prod_{2\leq k\leq N}A(k)(\phi_{k-1},\psi_{k-1})\right)\cdot\begin{pmatrix}
						F_{SU(N-1)}(\phi_{N},\ldots,\phi_{\frac{N(N-1)}{2}},\psi_{N},\ldots,\psi_{\frac{N(N-1)}{2}},\omega_1,\ldots,\omega_{N-2})&0\\
						0&1
					\end{pmatrix} e^{\lambda_{N^2-1}\omega_{N-1}}},
			\end{split}
		\end{equation}
		where $A(k)(x,y):=e^{\lambda_{k^2-1}x}e^{\lambda_{k^2-2}y}$. Here we denote the product as $$\prod_{2\leq k\leq N}A(k)(\phi_{k-1},\psi_{k-1}):=A(2)(\phi_1,\psi_1)\cdot \cdots\cdot A(N)(\phi_{N-1},\psi_{N-1}).$$ This mapping is a diffeomorphism onto its image, which is $SU(N)$ up to a measure zero set. In addition, when extending $F_{SU(N)}$ to the mapping $\tilde{F}_{SU(N)}:[0,\pi)^{\frac{N(N-1)}{2}}\times \left[0,\frac{\pi}{2}\right]^{\frac{N(N-1)}{2}}\times[0,2\pi)^{N-1}\rightarrow SU(N)$, the image of $\tilde{F}_{SU(N)}$ is $SU(N)$. Finally, the Haar measure is given by
		\begin{align*}
			dg_{SU(N)} = C_N\left(\prod_{j=1}^{N-1}2\sin^{2j-1}(\psi_j)\cos(\psi_j)\right) d\phi_1\ldots d\phi_{N-1}d\psi_1\ldots d\psi_{N-1} dg_{SU(N-1)}d\omega_N
		\end{align*}
		for some constant $C_N>0$.
	\end{lemma}
	\begin{proof}
		Let $G=SU(N)$, and $\mathfrak{g}=\mathfrak{su}(N)$. We apply Theorem \ref{thm:KAK_decomp} inductively. We start with $SU(2)$. Then the lemma is restating the ordinary Euler Angles, see for example \cite{Vilenkin}. Now let the parameterization be true for $SU(N-1)$, then we show it for $SU(N)$, with $N\geq 3$. We apply Theorem \ref{thm:KAK_decomp} to $SU(N)$, which is simply connected with finite center. For the involutive automorphism, we choose the inner automorphism $$\theta=\Ad\left[e^{\frac{\pi i}{N}}\begin{pmatrix}\mathbf{1}_{N-1}&\\
			&-1 
		\end{pmatrix}\right].$$ 
		With this choice of involution, we find 
		\begin{align*}
			\mathfrak{k}&=\mathrm{span}_\R(\lambda_1,\ldots,\lambda_{(N-1)^2-1},\lambda_{N^2-1})\simeq\mathfrak{su}(N-1)\oplus\mathfrak{u}(1),\\
			\mathfrak{p}&=\mathrm{span}_\R(\lambda_{(N-1)^2},\ldots,\lambda_{N^2-2}).
		\end{align*}  The corresponding connected subgroup $K=\langle \exp(\mathfrak{k})\rangle$ is then given by $K\simeq S(U(N-1)\times U(1))$ and can be embedded in $G$ as a submanifold of the form $$K=\left\{\left.\begin{pmatrix}
			A&\\
			&1
		\end{pmatrix} e^{\omega_{N-1}\lambda_{N^2-1}}\right|\omega_{N-1}\in\left[0,2\pi\right),A\in SU(N-1)\right\} \simeq SU(N-1)\times U(1).$$ Note that $K\simeq SU(N)\times U(1)$ is as manifolds, not as groups.
		We take as maximal abelian subalgebra $\mathfrak{a}=\R\lambda_{N^2-2}$. To find the roots, we consider the maximal abelian subalgebra $\mathfrak{h}$ in $\mathfrak{g}_\C=\mathfrak{sl}(N,\C)$ such that $\mathfrak{a}_\C\subseteq \mathfrak{h}$. We choose $$\mathfrak{h}=\bigoplus_{j=2}^{N-2}\C\lambda_{j^2-1}\oplus \C(2\lambda_{(N-1)^2-1}+\lambda_{N^2-1})\oplus \C\lambda_{N^2-2}.$$
		A direct calculation shows that 
		\begin{align*}
			\Delta_\mathfrak{p} =\{&2i\alpha,-2i\alpha,e_1+\beta\pm i\alpha,-e_1-\beta\pm i\alpha,\ldots,e_{N-3}+\beta\pm i\alpha,-e_{N-3}-\beta\pm i\alpha,\\
			&e_{N-2}-3\beta\pm i\alpha,-e_{N-2}+3\beta\pm i\alpha\}
		\end{align*}
		where the linear functional $e_i:M_N(\C)\rightarrow\C$ is given by $e_i(E_{jk})=\delta_{i,j}\delta_{j,k}$, where $E_{jk}$ is the matrix whose entries have zeroes everywhere except the $j,k$-th entry, which is 1. In addition, $$\alpha:M_N(\C)\rightarrow\C,\quad \alpha(E_{jk})=\delta_{j,N-1}\delta_{k,N}$$ and $$\beta:M_N(\C)\rightarrow\C,\quad \beta(E_{jk})=\delta_{j,N}\delta_{k,N}$$ which means in particular $\alpha(\lambda_{N^2-2})=1$ and $\beta(2\lambda_{(N-1)^2-1}+\lambda_{N^2-1})=-1$. We choose a set of positive roots such that the set $\Delta_\mathfrak{p}^+$ is given by
		$$\Delta_\mathfrak{p}^+ =\{2i\alpha,e_1+\beta\pm i\alpha,\ldots,e_{N-3}+\beta\pm i\alpha, e_{N-2}-3\beta\pm i\alpha,\}.$$
		Note that for any $\mu\in \Delta_\mathfrak{p}^+$ we get $\mu(H)=2i\alpha(H)$ or $\mu(H)=\pm i\alpha(H)$ for any $H\in\mathfrak{a}$. Thus we see that $$\mathcal{A}=\left\{c\lambda_{N^2-2}|c\in \left[0,\frac{\pi}{2}\right]\right\}.$$
		A small computation shows 
		\begin{align*}
			M=Z_K(\mathfrak{a})&=\left\{\left.\begin{pmatrix}
				B&&\\
				& e^{-ix}&\\
				&& e^{-ix}
			\end{pmatrix}\right| B\in U(N-2), x\in[0,2\pi) \text{ such that } \det(B)e^{-2ix}=1\right\}\\
			&\simeq S(U(N-2)\times U(1)).
		\end{align*}
		Applying all of this to Theorem \ref{thm:KAK_decomp} we get that the mapping 
		\begin{align*}
			(K/M)\times \exp(\mathcal{A})\times K \rightarrow SU(N), \qquad (kM,\exp(c\lambda_{N^2-2}),l)\mapsto \exp(\Ad_G(k)H)l = k\exp(c\lambda_{N^2-2})k^{-1}l
		\end{align*}
		is surjective, and a diffeomorphism up to a measure zero set if we replace $\mathcal{A}$ with $\mathrm{int}(\mathcal{A})$, which in this case is $\mathrm{int}(\mathcal{A})=\{c\lambda_{N^2-2}\,|\,c\in(0,\frac{\pi}{2})\}$. Applying the induction hypothesis to the element $k^{-1}l\in K$, we see that the mapping 
		\begin{align*}
			(K/M)\times \exp(\mathcal{A})\times [0,\pi)^{\frac{(N-1)(N-2)}{2}}\times \left[0,\frac{\pi}{2}\right]^{\frac{(N-1)(N-2)}{2}}\times [0,2\pi)^{N-2}\rightarrow SU(N)
		\end{align*}
		given by
		\begin{align}\label{eq:SU(N)_parametrization_the_map_that_actually_does_the_parametrization}
			\begin{split}
				(kM,\exp(\psi_{N-1}\lambda_{N^2-2}),\phi_{N},&\ldots,\phi_{\frac{N(N-1)}{2}},\psi_{N},\ldots,\psi_{\frac{N(N-1)}{2}},\omega_2,\ldots,\omega_{N-1})\\
				&\mapsto k\exp(\psi_{N-1}\lambda_{N^2-2})\begin{pmatrix}
					F_{SU(N-1)}(\phi_N,\ldots,\omega_{N-2})&0\\
					0&1
				\end{pmatrix}e^{\lambda_{N^2-1}\omega_{N-1}}\end{split}
		\end{align}
		is surjective, and if we replace $\mathcal{A}$ with $\mathrm{int}(\mathcal{A})$ and $\left[0,\frac{\pi}{2}\right]$ by $\left(0,\frac{\pi}{2}\right)$ it is a diffeomorphism upon its image, which is $SU(N)$ up to a measure zero set. We thus see that, getting a parametrization of $kM\in K/M$, i.e. a set of elements $k\in K$ such that the elements $kM$ uniquely describe the manifold $K/M$ in a smooth way, concludes the proof. The rest of the proof is thus dedicated to finding this specific subset of $k\in K$.
		
		To describe $K/M$, recall that by the induction hypothesis all elements $k\in K\simeq  SU(N-1)\times U(1)$ can be written as
		\begin{align*}
			\begin{split}
				k=&\left(\prod_{2\leq k\leq N-1}A(k)(\tilde{\phi}_{k-1},\tilde{\psi}_{k-1})\right)\cdot\\
				&\scalemath{0.9}{\begin{pmatrix}
						F_{SU(N-2)}(\tilde{\phi}_{N-1},\ldots,\tilde{\phi}_{\frac{(N-1)(N-2)}{2}},\tilde{\psi}_{N-1},\ldots,\tilde{\psi}_{\frac{(N-1)(N-2)}{2}},\tilde{\omega}_1,\ldots,\tilde{\omega}_{N-3})&0&0\\
						0&1&0\\
						0&0&1
					\end{pmatrix} e^{\lambda_{(N-1)^2-1}\tilde{\omega}_{N-2}}e^{\lambda_{N^2-1}\tilde{\omega}_{N-1}}}
			\end{split}	
		\end{align*}
		for some $\tilde{\phi}_i,\tilde{\psi}_j$ and $\tilde{\omega}_k$. Note that the matrix element $$m:=\begin{pmatrix}
			U&&\\
			&e^{i\frac{\tilde{\omega}_{N-2}}{2}}&\\
			&&e^{i\frac{\tilde{\omega}_{N-2}}{2}}
		\end{pmatrix}$$ lies in $M$, where $$U = \begin{pmatrix}
			\mathbf{1}_{N-3}&\\
			&e^{-i\tilde{\omega}_{N-2}}
		\end{pmatrix}[F_{SU(N-2)}(\tilde{\phi}_{N-1}\ldots,\tilde{\omega}_{N-3})]^{-1}\in U(N-2).$$
		This shows that $$gm = \left(\prod_{2\leq k\leq N-1}A(k)(\tilde{\phi}_{k-1},\tilde{\psi}_{k-1})\right)e^{\lambda_{N^2-1}(\tilde{\omega}_{N-1}-\frac{\tilde{\omega}_{N-2}}{2})}.$$ 
		In other words, we see that the set $$X:=\left\{\left.\left(\prod_{2\leq k\leq N-1}A(k)(\phi_{k-1},\psi_{k-1})\right)e^{\lambda_{N^2-1}\phi_{N-1}}\right|\phi_1,\ldots,\phi_{N-1}\in[0,\pi)\text{ and }\psi_1,\ldots,\psi_{N-2}\in\left[0,\frac{\pi}{2}\right]\right\}$$ is a  candidate for parametrizing $K/M$ up to a measure zero set, i.e. $K/M = \{gM\,|\,g\in X\}$ up to a measure zero set. To prove that it is a parametrization, let $g,h\in X$. We show that $gM\cap hM=\emptyset.$ In other words, if there exists $m\in M$ such that $gm=h$, then $g=h$. We restrict ourselves to the case $N=3$, for the higher dimensional cases can be reduced to the case $N=3$ by considering the lower-right $3\times 3$ matrix in the $gm=h$ equation. Let $g\in X$ be parametrized by $\phi_1,\phi_2\in [0,\pi)$ and $\psi\in [0,\pi/2]$, and $h\in X$ by $\tilde{\phi}_1,\tilde{\phi}_2\in[0,\pi)$ and $\tilde{\psi}\in[0,\pi/2]$. Let $m\in M$ which in this case can be written as $m=\mathrm{diag}(e^{2ix},e^{-ix},e^{-ix})$ for $x\in[0,2\pi)$. Then we get the equation given by
		\begin{align*}
			\begin{pmatrix}
				e^{i\phi_{1}}&&\\
				&e^{-i\phi_{1}}&\\
				&&1
			\end{pmatrix}
			\begin{pmatrix}
				\cos(\psi_{1})&\sin(\psi_1)&\\
				-\sin(\psi_1)&\cos(\psi_1)&\\
				&&1
			\end{pmatrix}
			\begin{pmatrix}
				1&&\\
				&e^{i\phi_2}&\\
				&&e^{-i\phi_2}
			\end{pmatrix}
			\begin{pmatrix}
				e^{2ix}&&\\
				&e^{-ix}&\\
				&&e^{-ix}
			\end{pmatrix} = \\
			\begin{pmatrix}
				e^{i\tilde{\phi}_{1}}&&\\
				&e^{-i\tilde{\phi}_{1}}&\\
				&&1
			\end{pmatrix}
			\begin{pmatrix}
				\cos(\tilde{\psi}_{1})&\sin(\tilde{\psi}_1)&\\
				-\sin(\tilde{\psi}_1)&\cos(\tilde{\psi}_1)&\\
				&&1
			\end{pmatrix}
			\begin{pmatrix}
				1&&\\
				&e^{i\tilde{\phi}_2}&\\
				&&e^{-i\tilde{\phi}_2}
			\end{pmatrix}.
		\end{align*}
		The equation in the lower right component gives $e^{-i(\phi_2+x)}=e^{-i\tilde{\phi}_2}$, so $\phi_2-\tilde{\phi}_2=x+2\pi k$ for some $k\in\Z$. Putting this in gives
		\begin{align*}
			\begin{pmatrix}
				e^{i(\phi_{1}-\tilde{\phi}_1)+2ix}\cos(\psi_{1})&e^{i(\phi_{1}-\tilde{\phi}_1)-2ix}\sin(\psi_1)&\\
				-e^{-i(\phi_{1}-\tilde{\phi}_1)+2ix}\sin(\psi_1)&e^{-i(\phi_{1}-\tilde{\phi}_1)-2ix}\cos(\psi_1)&\\
				&&1
			\end{pmatrix} =
			\begin{pmatrix}
				\cos(\tilde{\psi}_{1})&\sin(\tilde{\psi}_1)&\\
				-\sin(\tilde{\psi}_1)&\cos(\tilde{\psi}_1)&\\
				&&1
			\end{pmatrix}.
		\end{align*}
		Note that the right-hand side is a real matrix. Hence all exponentials should be either $1$ or $-1$. Now $\psi_1,\tilde{\psi}_1\in[0,\pi/2]$, so the sine and cosine are both non-negative and injective on this interval. Therefore we must have $\phi_1-\tilde{\phi}_1+2x=2\pi l$ and $\phi_1-\tilde{\phi}_1-2x=2\pi l'$ with $l,l'\in \mathbb{Z}$. In other words, $x=\frac{(l-l')\pi}{2}$ and $\phi_1-\tilde{\phi}_1 = (l+l')\pi$. Now since $\phi_1,\phi_1'\in [0,\pi)$ we have that $l=-l'$, hence $\phi_1=\tilde{\phi}_1$. This also means that $x=\pi l$. But now $\phi_2-\tilde{\phi}_2= \pi l +2\pi k$, and remember that $\phi_2,\tilde{\phi}_2\in [0,\pi)$ which means that that can only be true if $l=k=0$. In other words, we have $g=h$. This shows that $X$ is in bijection with $K/M$, and thus by replacing $hM$ with $h\in X$ one gets the surjectivity of the map $F_{SU(N)}$ as map $[0,\pi)^{\frac{N(N-1)}{2}}\times \left[0,\frac{\pi}{2}\right]^{\frac{N(N-1)}{2}}\times[0,2\pi)^{N-1}\rightarrow SU(N)$. 
		
		To show that $F_{SU(N)}$ as map $[0,\pi)^{\frac{N(N-1)}{2}}\times \left(0,\frac{\pi}{2}\right)^{\frac{N(N-1)}{2}}\times[0,2\pi)^{N-1}\rightarrow SU(N)$ is a diffeomorphism upon its image, we note that, by previous arguments, it is enough to show that the map $f:Y\rightarrow K/M$ given by $f(g)= gM$ is a diffeomorphism unto its image, where $$Y:=\left\{\left.\left(\prod_{2\leq k\leq N-1}A(k)(\phi_{k-1},\psi_{k-1})\right)e^{\lambda_{N^2-1}\phi_{N-1}}\right|\phi_1,\ldots,\phi_{N-1}\in[0,\pi)\text{ and }\psi_1,\ldots,\psi_{N-2}\in\left(0,\frac{\pi}{2}\right)\right\}.$$ It is clear that $f$ is smooth if we endow $Y\subseteq K$ with the subset topology. An extensive but straightforward calculation shows the tangent map $T_xf:T_xY\rightarrow T_{xM}(K/M)$ given by $T_el_x (H)\mapsto T_{eM}\tau_x(H+\mathfrak{m})$ is surjective, where $\tau_x$ is the diffeomorphism $\tau_x:K/M\rightarrow K/M$ given by $\tau_x(gM)=xgM$ and $\mathfrak{m}=\mathrm{Lie}(M)$. Since $$\dim T_gY = 2N-3 = T_{gM}(K/M)$$ we see that $T_xf$ is bijective, so $f$ is in fact a diffeomorphism upon its image, which is $K/M$ up to a measure zero set. This proves Equation (\ref{eq:SU(N)_parametrization_new_version}).
		
		To show the form of the Jacobian, we note that by Theorem \ref{thm:KAK_decomp} we have 
		\begin{align*}
			dg=|J(\exp(H))|dHdk_Mdk &= \left|\prod_{\alpha\in\Delta_\mathfrak{p}^+}\sin(\alpha(\psi_{N-1}\lambda_{N^2-2}))\right|d\psi_{N-1}dk_Mdk\\
			&=\left|\sin(2\psi_{N-1})\prod_{j=1}^{N-2}\left(\sin(\psi_{N-1})\sin(-\psi_{N-1})\right)\right|d\psi_{N-1}dk_Mdk\\
			&=2\sin^{2(N-1)-1}(\psi_{N-1})\cos(\psi_{N-1})d\psi_{N-1}dk_Mdk.
		\end{align*}
		In addition, since $K\simeq SU(N-1)\times U(1)$, the Haar measure on $K$ decomposes as $dk=dg_{SU(N-1)}d\omega_{N-1}$. The decomposition of the measure $dk_M$ proceeds in the same way as in the original proof of \cite{Zwart}, and thus $$dk_M=\left(\prod_{j=1}^{N-2}2\sin^{2j-1}(\psi_j)\cos(\psi_{j})\right)d\phi_1 d\phi_2\ldots d\phi_{N-1}d\psi_1\ldots d\psi_{N-2}$$ up to a constant, proving the lemma. 
	\end{proof}
	
	As in our previous work, we are interested in the finite-type functions of $SU(N)$. We recall:
	\begin{definition}
		Let $G$ be a compact Lie group. A function $f:G\rightarrow \C$ is called a \emph{finite-type function} if it can be written as a finite linear combination of matrix coefficients of irreducible representations, i.e. $$f(x)=\sum_{j=1}^n \mathrm{Tr}(a_j\pi_j(x))$$ where $(\pi_j,V_j)$ is an irreducible representation of $G$, and $a_j\in \mathrm{End}(V_j)$.
	\end{definition}
	\begin{theorem}\cite[Thm. 8.2.3]{Procesi}\label{thm:matrix_coefficients_generated_by_matrix_entries}
		Let $G\subseteq U(N)$ be a connected compact Lie group. Then the ring of finite-type functions on $G$ is generated by the matrix entries and the inverse of the determinant.
	\end{theorem}
	
	With the parametrization in Lemma \ref{lemma:Euler_Angles_SU(N)_revisited}, it is clear that the finite-type functions, as noted in Equation (2.3) of \cite{Zwart}, are the same as the ones we would get from this parametrization. However, Lemma 2.7 in \cite{Zwart} can be improved.
	
	\begin{definition}\label{def:g_{SU(N)}}
		Let $G=SU(N)$ and let $g\in G$ be such that there exist parameters $\phi_1,\ldots,\omega_{N-1}$ such that $F_{N}(\phi_1,\ldots,\omega_{N-1})=g$. By Lemma \ref{lemma:Euler_Angles_SU(N)_revisited} this is true for almost all $g\in G$. For these $g\in G$ we will use the shorthand notation $g_{SU(1)}:=1$ and  $$g_{SU(n)}:=F_{SU(n)}(\phi_{\frac{N(N-1)}{2}-(\frac{n(n-1)}{2}-1)},\ldots,\phi_{\frac{N(N-1)}{2}},\psi_{\frac{N(N-1)}{2}-(\frac{n(n-1)}{2}-1)},\ldots,\psi_{\frac{N(N-1)}{2}}, \omega_{1},\ldots,\omega_{n-1})$$ for $2\leq n\leq N$. Note that $g_{SU(n)}\in SU(n)$.
	\end{definition}
	
	\begin{lemma}\label{lemma:rewrite_integral_for_SU(N)}
		Let $N\geq2$, let $g\in SU(N)$ and let $N_-=\frac{N(N-1)}{2}$. Define the finite-type function $$f^{SU(2)}(g_{SU(2)}):=c e^{ik_{N_-}\phi_{N_-}}\sin^{m_{N_-}}(\psi_{N_-})\cos^{n_{N_-}}(\psi_{N_-})e^{il_{1}\omega_{1}},$$ for some $c\in\C$, $m_{N_-}, n_{N_-}\in \N_0$ and $k_{N_-},l_1\in \Z$, and define the finite-type function $f^{SU(N)}$ recursively as
		\begin{align*}
			f^{SU(N)}(g_{SU(N)})=&e^{ik_1\phi_1}\sin^{m_1}(\psi_1)\cos^{n_1}(\psi_1)\cdots e^{ik_{N-1}\phi_{N-1}}\sin^{m_{N-1}}(\psi_{N-1})\cos^{n_{N-1}}(\psi_{N-1})\cdot\\
			&f^{SU(N-1)}(g_{SU(N-1)})e^{il_{N-1}\omega_{N-1}}
		\end{align*} where $k_1,\ldots,k_{N-1},l_{N-1}\in\Z$ and $m_1,\ldots,m_{N-1},n_{1},\ldots,n_{N-1}\in\N_0$. Then 
		\begin{align}\label{eq:Integral_evaluated_SU(N)_finite_type}
			\begin{split}\int_{SU(N)}& f^{SU(N)}(g)dg= 2\pi^Nc\delta_{k_1,0}\ldots\delta_{k_{N-1},0}\delta_{l_{N-1},0}\int_{SU(N-1)}f^{SU(N-1)}(g_{SU(N-1)}) dg_{SU(N-1)}\cdot\\
				&\int_{[0,1]^{N-1}}x_1^{m_1}(1-x_1^2)^{\frac{n_1}{2}}\cdots x_{N-1}^{m_{N-1}}(1-x_{N-1}^2)^{\frac{n_{N-1}}{2}} \tilde{J}_{SU(N)}(x_1,\ldots,x_{N-1}) dx_1\ldots dx_{N-1}.\end{split}
		\end{align} Here $dg_{SU(N-1)}$ is the Haar measure on $SU(N-1)$, and $\tilde{J}_{SU(N)}$ is given by
		\begin{align*}
			\tilde{J}_{SU(N)}(x_1,\ldots,x_N)=2^{N-1}C_N \prod_{j=1}^{N-1}x_j^{2j-1}.
		\end{align*}
		where $C_N$ is some constant.
	\end{lemma}
	\begin{remark}
		The integral over the $x$-variables in Equation $(\ref{eq:Integral_evaluated_SU(N)_finite_type})$ can be evaluated by noting that the integral can be split into multiple one-dimensional integrals, i.e. 
		\begin{align*}
			\int_{[0,1]^{N-1}}x_1^{m_1}(1-x_1^2)^{\frac{n_1}{2}}\cdots x_{N-1}^{m_{N-1}}&(1-x_{N-1}^2)^{\frac{n_{N-1}}{2}} \tilde{J}_{SU(N)}(x_1,\ldots,x_{N-1}) dx_1\ldots dx_{N-1}\\
			&=2^{N-1}C_N\prod_{j=1}^{N-1}\int_0^1 x_j^{m_j+2j-1}(1-x_j^2)^{\frac{n_j}{2}}\,dx_j.
		\end{align*}
		Now we can calculate the latter integrals, by substituting $t=x_j^2$ and noticing the definition of the Beta function, which can be expressed as a quotient of Gamma functions, giving
		\begin{align*}
			2^{N-1}C_N\prod_{j=1}^{N-1}\int_0^1 x_j^{m_j+2j-1}(1-x_j^2)^{\frac{n_j}{2}}\,dx_j &= C_N\prod_{j=1}^{N-1}\int_0^1 t^{\frac{m_j}{2}+j-1}(1-t)^{\frac{n_j}{2}}\,dt\\
			&= C_N\prod_{j=1}^{N-1}\frac{\Gamma(\frac{m_j}{2}+j)\Gamma(\frac{n_j}{2}+1)}{\Gamma(\frac{m_j+n_j}{2}+j+1)}.
		\end{align*}
		Albeit useful, we will not pursue the actual evaluation of these integrals in this paper, for Lemma \ref{lemma:rewrite_integral_for_SU(N)} is enough for us to produce the desired results.
	\end{remark}

	\begin{proof}
		We note that our definition of $f^{SU(N)}$ covers all monomials in the ring of finite-type functions by Theorem \ref{thm:matrix_coefficients_generated_by_matrix_entries}. Now to show the equality, we make extensive use of the properties of the Haar measure and Lemma \ref{lemma:Euler_Angles_SU(N)_revisited}. We remind ourselves that $G=SU(N)$ is compact, hence the Haar measure is unimodular, i.e. $$\int_{G}f(gy)dg=\int_{G}f(g)dg=\int_{G}f(yg)dg$$ for any $y\in G$. This must restrict the possible values of the integral. The idea of the proof is then to choose specific $y\in G$ in such a way that the result follows. Let $g_{SU(N)}\in G$ be as in Definition \ref{def:g_{SU(N)}}. Choosing $y= e^{t\lambda_3}$ and describing the element $e^{t\lambda_3}g_{SU(N)}$ using Lemma \ref{lemma:Euler_Angles_SU(N)_revisited}, we see that 
		$$e^{t\lambda_3}g_{SU(N)}=F_{SU(N)}(\phi_1+t,\phi_2,\ldots,\phi_{\frac{N(N-1)}{2}},\psi_1,\ldots,\psi_{\frac{N(N-1)}{2}},\omega_1,\ldots,\omega_{N-1}).$$
		Hence in the integral it translates to	
		$$\int_G f(g) dg = \int_G f(e^{t\lambda_3}g) dg = e^{ik_1t}\int_G f(g)dg.$$
		This is true for all $t\in\R$, thus we must have $$k_1=0.$$
		In a similar fashion, considering $y= e^{t\lambda_{N^2-1}}$, we find $$g_{SU(N)}e^{t\lambda_{N^2-1}}=F_{SU(N)}(\phi_1,\ldots,\omega_{N-2},\omega_{N-1}+t),$$ and thus $$\int_G f(g)\, dg = \int_G f(ge^{t\lambda_{N^2-1}})\, dg = e^{in_{N-1}t}\int_G f(g) \,dg,$$ which can only be true if $$n_{N-1}=0.$$
		To get a similar result for the other parameters, more extensive computations are needed. Note that the following equation holds
		\begin{align}\label{eq:Integral_SU(N)_giving_deltas_pulling_exponentials_through_1}
			\Ad\left(\begin{pmatrix}
				e^{it}&0\\
				0&1
			\end{pmatrix}\right)\begin{pmatrix}
				\cos(\psi)&\sin(\psi)\\
				-\sin(\psi)&\cos(\psi)
			\end{pmatrix} = \Ad\left(\begin{pmatrix}
				e^{it/2}&0\\0&e^{-it/2}
			\end{pmatrix}\right)\begin{pmatrix}
				\cos(\psi)&\sin(\psi)\\-\sin(\psi)&\cos(\psi)
			\end{pmatrix},
		\end{align}
		and similarly
		\begin{align}\label{eq:Integral_SU(N)_giving_deltas_pulling_exponentials_through_2}
			\Ad\left(\begin{pmatrix}
				1&0\\
				0&e^{it}
			\end{pmatrix}\right)\begin{pmatrix}
				\cos(\psi)&\sin(\psi)\\
				-\sin(\psi)&\cos(\psi)
			\end{pmatrix} = \Ad\left(\begin{pmatrix}
				e^{-it/2}&0\\0&e^{it/2}
			\end{pmatrix}\right)\begin{pmatrix}
				\cos(\psi)&\sin(\psi)\\-\sin(\psi)&\cos(\psi)
			\end{pmatrix}.
		\end{align}
		With these equalities, we see for example that
		\begin{align*}
			e^{t\lambda_8}e^{\psi_1\lambda_2} = e^{-\frac{t}{2}\lambda_3}e^{\psi_1\lambda_2}e^{\frac{t}{2}\lambda_3}e^{t\lambda_8}
		\end{align*}
		and similarly
		\begin{align}\label{eq:Addendum_proof_rewriting_integral_pulling_exps}
			\begin{split}e^{t\lambda_8}e^{\psi_1\lambda_2}e^{\phi_2\lambda_{8}}e^{\psi_2\lambda_{7}} &= e^{-\frac{t}{2}\lambda_3}e^{\psi_1\lambda_2}e^{\frac{t}{2}\lambda_3}e^{(\phi_2+t)\lambda_{8}}e^{\psi_2\lambda_{7}}\\
				&= e^{-\frac{t}{2}\lambda_3}e^{\psi_1\lambda_2}e^{(\phi_2+t)\lambda_{8}}e^{\frac{t}{2}\lambda_3}e^{\psi_2\lambda_{7}}\\
				&= e^{-\frac{t}{2}\lambda_3}e^{\psi_1\lambda_2}e^{(\phi_2+t)\lambda_{8}}e^{-\frac{t}{4}\lambda_8}e^{\psi_2\lambda_{7}}e^{\frac{t}{4}\lambda_8}e^{\frac{t}{2}\lambda_3}\\
				&=e^{-\frac{t}{2}\lambda_3}e^{\psi_1\lambda_2}e^{(\phi_2+\frac{3t}{4})\lambda_{8}}e^{\psi_2\lambda_{7}}e^{\frac{t}{4}\lambda_8}e^{\frac{t}{2}\lambda_3}.\end{split}
		\end{align}
		The strategy here is that we are pulling $e^{t\lambda_8}$ to the right, using the commutation relations given in Equation $(\ref{eq:Integral_SU(N)_giving_deltas_pulling_exponentials_through_1})$ and $(\ref{eq:Integral_SU(N)_giving_deltas_pulling_exponentials_through_2})$, until we find an exponential of the form $e^{\phi\lambda_8}$ for some $\phi$. Then we will try to pull the last new exponential we got from the latest commutation relation to the right (e.g. in Equation (\ref{eq:Addendum_proof_rewriting_integral_pulling_exps}) we mean $e^{\frac{t}{2}\lambda_3}$), using the same commutation relations again. This process will continue until we have pulled every matrix this way all the way to the right.
		
		Now consider $e^{t\lambda_8}g$ for almost all $g\in G$. Then using Equation (\ref{eq:Addendum_proof_rewriting_integral_pulling_exps}), we see that
		\begin{align*}
			e^{t\lambda_{8}}g_{SU(N)}&= e^{t\lambda_8}F(\phi_1,\ldots,\omega_{N-1})\\
			&=e^{t\lambda_{8}}e^{\phi_1\lambda_{3}}e^{\psi_1\lambda_{2}}e^{\phi_2\lambda_{8}}e^{\psi_2\lambda_{7}}\ldots e^{\phi_{N-1}\lambda_{N^2-1}}e^{\psi_{N-1}\lambda_{N^2-2}} g_{SU(N-1)} e^{\omega_{N-1}\lambda_{N^2-1}}\\
			&=e^{(\phi_1-\frac{t}{2})\lambda_{3}}e^{\psi_1\lambda_{2}}e^{(\phi_2+t)\lambda_{8}}e^{-\frac{t}{4}\lambda_{8}}e^{\psi_2\lambda_{7}}e^{\frac{t}{4}\lambda_8}e^{\frac{t}{2}\lambda_{3}}\ldots e^{\phi_{N-1}\lambda_{N^2-1}}e^{\psi_{N-1}\lambda_{N^2-2}} g_{SU(N-1)} e^{\omega_{N-1}\lambda_{N^2-1}}\\
			&=e^{(\phi_1-\frac{t}{2})\lambda_{3}}e^{\psi_1\lambda_{2}}e^{(\phi_2+\frac{3t}{4})\lambda_{8}}e^{\psi_2\lambda_{7}}e^{\frac{t}{4}\lambda_8}\ldots e^{\phi_{N-1}\lambda_{N^2-1}}e^{\psi_{N-1}\lambda_{N^2-2}} e^{\frac{t}{2}\lambda_{3}}g_{SU(N-1)} e^{\omega_{N-1}\lambda_{N^2-1}}.
		\end{align*}
		where in the last equation we used that $e^{\frac{t}{2}\lambda_3}$ commutes with all elements to the right up and till $g_{SU(N-1)}$. Pushing $e^{\frac{t}{4}\lambda_8}$ to the right, using the above mentioned strategy we find that 
		$$e^{t\lambda_8}g_{SU(N)}=e^{\tilde{\phi}_1\lambda_3}e^{\psi_1\lambda_{2}}e^{\tilde{\phi}_2\lambda_{8}}e^{\psi_2\lambda_{7}}\ldots e^{\tilde{\phi}_{N-1}\lambda_{N^2-1}}e^{\psi_{N-1}\lambda_{N^2-2}} kg_{SU(N-1)} e^{\omega_{N-1}\lambda_{N^2-1}}$$
		where $\tilde{\phi}_1=\phi_1-\frac{t}{2}$, $\tilde{\phi}_2=\phi_2+\frac{3}{4}t$ and $\tilde{\phi}_j=\phi_j-\frac{t}{2^j}$ for $3\leq j\leq N-1$ and $$k:=\begin{pmatrix}
			e^{it/2}&&&&\\
			&e^{-it/4}&&&\\
			&&\ddots&&\\
			&&&e^{-it/2^{N-1}}&\\
			&&&&e^{-it/2^{N-1}}
		\end{pmatrix}.$$  
		Note $k\in K$. This way, we get
		\begin{align*}
			\int_G f(g)dg = e^{it\left(-\frac{k_1}{2}+\frac{3k_2}{4}-\sum_{j=3}^{N-1}\frac{k_j}{2^j}\right)}\int_G& ce^{ik_1\phi_1}\sin^{m_1}(\psi_1)\cos^{n_1}(\psi_1)\cdots e^{ik_{N-1}\phi_{N-1}}\cdot\\
			&\sin^{m_{N-1}}(\psi_{N-1})\cos^{n_{N-1}}(\psi_{N-1})f^{SU(N-1)}(kg_{SU(N-1)})e^{il_N\omega_{N-1}}dg.
		\end{align*}
		We recall that the Haar measure $dg$ decomposes into $$dg= J(\psi_1,\ldots,\psi_{N-1}) d\phi_1\ldots d\phi_{N-1}d\psi_1\ldots d\psi_{N-1} dk d\omega_{N-1}$$ by Lemma \ref{lemma:Euler_Angles_SU(N)_revisited}, where $J:(0,\frac{\pi}{2})^{N-1}\rightarrow \C$ is given by $J(\psi_1,\ldots,\psi_{N-1}):=\prod_{j=1}^{N-1}2\sin^{2j-1}(\psi_j)\cos(\psi_j)$. We see that $dg$ decomposes as the Haar measure on $K$ times other measures, so $$\int_{SU(N-1)}f^{SU(N-1)}(kg_{SU(N-1)})dg_{SU(N-1)}=\int_{SU(N-1)}f^{SU(N-1)}(g_{SU(N-1)})dg_{SU(N-1)},$$ giving
		\begin{align*}
			\int_G f(g)dg = e^{it\left(-\frac{k_1}{2}+\frac{3k_2}{4}-\sum_{j=3}^{N-1}\frac{k_j}{2^j}\right)}\int_G f(g) dg.
		\end{align*}
		Since this is true for all $t\in \R$, we can conclude that $$-\frac{k_1}{2}+\frac{3k_2}{4}-\sum_{j=3}^{N-1}\frac{k_j}{2^j}=0.$$ The same procedure can be performed by considering $e^{t\lambda_{i^2-1}}g$ instead of $e^{t\lambda_8}g$ for $i=3,\ldots, N-1$. Applying the same steps, we get the following set of equations
		\begin{align*}
			-\frac{k_{i-1}}{2}+\frac{3k_i}{4}-\sum_{j=i+1}^{N-1}\frac{k_j}{2^j}&=0 \qquad \forall i=2,\ldots,N-1,\\
			-\frac{k_{N-2}}{2}+\frac{3k_{N-1}}{4}&=0.
		\end{align*}
		We can collect the set of equations we have found in a linear system of the form $A\vec{k}=0$ where $\vec{k}=(k_1,\ldots,k_{N-1})^T$, and $$A:=\begin{pmatrix}
			1&0&0&0&0&\ldots&0\\
			-\frac{1}{2}&\frac{3}{4}&-\frac{1}{8}&-\frac{1}{16}&-\frac{1}{32}&\ldots&-\frac{1}{2^{N-1}}\\
			0&-\frac{1}{2}&\frac{3}{4}&-\frac{1}{8}&-\frac{1}{16}&\ldots&-\frac{1}{2^{N-2}}\\
			0&0&-\frac{1}{2}&\frac{3}{4}&-\frac{1}{8}&\ldots&-\frac{1}{2^{N-3}}\\
			\vdots&&\ddots&\ddots&\ddots&&\vdots\\
			0&\ldots&\ldots&0&-\frac{1}{2}&\frac{3}{4}&-\frac{1}{8}\\
			0&0&0&0&0&-\frac{1}{2}&\frac{3}{4}
		\end{pmatrix}.$$
		Using a recursive method by repeatedly developing to the left column, we can show that this Hessenberg matrix has determinant $\det(A)=\frac{N}{2^{N-1}}$ for $N\geq 4$, and thus is invertible. Thus $\vec{k}=0$. In other words, $$k_i=0\quad\forall i=1,\ldots, N-1.$$
		This way, the integral becomes	
		\begin{align*}
			\int_G f^{SU(N)}(g) =& \delta_{k_1,0}\ldots\delta_{k_{N-1},0}\delta_{n_{N-1},0}\int_{[0,\pi]^{N-1}}\int_{[0,2\pi]}\int_{[0,\frac{\pi}{2}]^{N-1}}\int_{SU(N-1)}c\sin^{m_1}(\psi_1)\cos^{n_1}(\psi_1)\cdots\\
			&\sin^{m_{N-1}}(\psi_{N-1})\cos^{n_{N-1}}(\psi_{N-1})f^{SU(N-1)}(g_{SU(N-1)}) J(\psi_1\,\ldots,\psi_{N-1})\, dg_{SU(N-1)}d\psi_1\ldots\\
			& d\psi_{N-1} d\omega_{N-1}d\phi_1\ldots d\phi_{N-1}\\
			=&2\pi^Nc\delta_{k_1,0}\ldots\delta_{k_{N-1},0}\delta_{n_{N-1},0}\left(\int_{SU(N-1)}f^{SU(N-1)}(g_{SU(N-1)})dg_{SU(N-1)}\right)\cdot\\
			&\int_{[0,\frac{\pi}{2}]^{N-1}}\sin^{m_1}(\psi_1)\cdots\cos^{n_{N-1}}(\psi_{N-1}) J(\psi_1\,\ldots,\psi_{N-1})\,d\psi_1\ldots d\psi_{N-1}.
		\end{align*}
		To complete the proof, we make use of the following equality
		\begin{align*}
			\int_{0}^{\pi/2}\sin^{k+p}(\psi)\cos^{l+q}(\psi)\;d\psi &=\int_0^1 x^{k+p}(1-x^2)^{\frac{l+q-1}{2}}\; dx.
		\end{align*}
		for any $k,p,q\in\N_0$ and $l\in\N$.
	\end{proof}
	Using the fact that $\delta_{k,0}=\frac{1}{2\pi}\int_{0}^{2\pi}e^{ik\phi}d\phi = \frac{1}{2\pi i}\int_{\mathbb{T}}z^k\frac{dz}{z},$ we immediately find the following corollary:
	\begin{corollary}\label{cor:rewrite_integral_for_SU(N)}
		Let $f^{SU(N)}$ be a finite-type function as in \cite[Equation (2.3)]{Zwart}, let $N\geq 2$, and denote $N_{\pm}=\frac{N(N\pm 1)}{2}$. Then 
		\begin{align}
			\begin{split}
				\int_{SU(N)} f^{SU(N)}(g)dg = C&\int_{[0,1]^{N_-}}\int_{\mathbb{T}^{N_+-1}}\widetilde{f^{SU(N)}}(x_1,\ldots,z_{N_+-1})\cdot\\ &\tilde{J}_{SU(N)}(x_1,\ldots,x_{N_-})\frac{dz_1}{z_1}\ldots \frac{dz_{N_+-1}}{z_{N_+-1}}dx_1\ldots dx_{N_-}.\end{split}
		\end{align} Here $C$ is some number that is independent of $f^{SU(N)}$, and $\widetilde{f^{SU(N)}}$ is defined recursively by $\widetilde{f^{SU(1)}}\equiv 1$ and 
		\begin{align}\label{eq:rewrite_integral_for_SU(N)_finite_type_function}
			\begin{split}\widetilde{f^{SU(N)}}(x_1,\ldots,z_{N_+}) := \sum_{j=1}^M\sum_{i=1}^Q& c_{ij}z_1^{k_{ij}^1}x_1^{m_{ij}^1}(1-x_1^2)^{\frac{n_{ij}^1}{2}}\cdots z_{N-1}^{k_{ij}^{N-1}}x_{N-1}^{m_{ij}^{N-1}}(1-x_{N-1}^2)^{\frac{n_{ij}^{N-1}}{2}}\\
				&\widetilde{f^{SU(N-1)}}(x_N,\ldots,x_{N_-},z_{N},\ldots,z_{N_+-2})\, z_{N_+-1}^{l_{ij}^{N-1}}.\end{split}
		\end{align} 
		Similarly $\tilde{J}_{SU(N)}$ is defined recursively as $\tilde{J}_{SU(1)}\equiv 1$ and 
		\begin{align*}
			\tilde{J}_{SU(N)}(x_1,\ldots,x_{N_-}):= 2^{N-1}C_N \left(\prod_{j=1}^{N-1}x_j^{2j-1}\right)\tilde{J}_{SU(N-1)}(x_N,\ldots,x_{N_-})
		\end{align*}
	\end{corollary}
	
	\begin{remark}
		We note that this lemma is an improvement to Lemma 2.7 in \cite{Zwart} for there are no roots in the $z$ variables here. This allows us to use $\mathbb{T}^n$ instead of $(S^1-\{1\})^n$, opening up for more tools to analyse these functions. 
	\end{remark}
	
	The absence of roots in the $z$-variables also allows us to redefine what an admissible function is:
	
	\begin{definition}\label{def:admissible_function}
		Let $k,l\in\N$ and $f: [0,1]^{k}\times \mathbb{T}^{l}\rightarrow\C$. We say $f$ is an \emph{admissible function} if $f$ can be written as $$f(x_1,\ldots,x_k,z_1,\ldots,z_{l})=\sum_{\vec{m}}c_{\vec{m}}(x)z^{\vec{m}},$$ where $\vec{m}=(m_1,\ldots,m_l)$ is a multi-index where $m_i\in\mathbb{Z}$, and $c_{\vec{m}}(x)\in \C[x_1,\sqrt{1-x_1^2},\ldots,x_k,\sqrt{1-x_k^2}]$ is a complex polynomial in $x_i$ and $\sqrt{1-x_i^2}$. We call the collection of $\vec{m}$ for which $c_{\vec{m}}\neq 0$ \emph{the spectrum of} $f$, and it will be denoted by $\mathrm{Sp}(f)$.
	\end{definition} 
	
	As expected in view of \cite{Zwart,Zwart2}, we have the following conjecture and Theorem, which is proven by using Lemma \ref{lemma:rewrite_integral_for_SU(N)} extensively and the last part of the proof of Theorem 2.11 in \cite{Zwart}. 
	
	\begin{conjecture}\label{con:xz-conjecture}
		Let $f:[0,1]^{\frac{N(N-1)}{2}}\times \mathbb{T}^{\frac{N(N+1)}{2}-1}\rightarrow\C$ be an admissible function. If $$\int_{[0,1]^{\frac{N(N-1)}{2}}}\int_{\mathbb{T}^{\frac{N(N+1)}{2}-1}}f^P \tilde{J}_{SU(N)} \,\,\frac{dz_1}{z_1}\ldots \frac{dz_{\frac{N(N+1)}{2}-1}}{z_{\frac{N(N+1)}{2}-1}}dx_1\ldots dx_{\frac{N(N-1)}{2}} = 0$$ for all $P\in \N$, then $\vec{0}$ does not lie in the convex hull of $\mathrm{Sp}(f)$.
	\end{conjecture}
	
	\begin{theorem}\label{thm:Mathieu_proven_assuming_XZ-conjecture_SU(N)}
		Assume Conjecture \ref{con:xz-conjecture} is true. Then the Mathieu Conjecture is true for $SU(N)$.
	\end{theorem}
	
	\section{Concerning the groups $Sp(N)$ and $G_2$}\label{sec:Applying_to_Sp(N)_and_G_2}
	In \cite{Zwart2}, we used the decomposition of $SU(N)$ to obtain results about $Sp(N)$ and $G_2$. Using Lemma \ref{lemma:Euler_Angles_SU(N)_revisited} and the techniques given in this paper instead of Lemma 2.5 in \cite{Zwart2}, we arrive at the following two conjectures and theorems: 
	
	\begin{conjecture}\label{con:xz-conjecture_Sp(N)}
		Let $f:[0,1]^{N^2}\times \mathbb{T}^{N(N+1)}\rightarrow\C$ be an admissible function in the sense of Definition \ref{def:admissible_function}. If $$\int_{[0,1]^{N(N-1)}}\int_{\mathbb{T}^{N(N+1)}}\int_0^1\int_0^{\xi_N}\ldots\int_0^{\xi_2}f^P \tilde{J}_{Sp(N)} \,\,d\xi_1\ldots d\xi_N\frac{dz_1}{z_1}\ldots \frac{dz_{N(N+1)}}{z_{N(N+1)}}dx_1\ldots dx_{N(N-1)} = 0$$ for all $P\in \N$, where \begin{align*}
			\tilde{J}_{Sp(N)}(x_1,\ldots,x_{N(N-1)},\xi_1,\ldots,\xi_N):=\,&\tilde{J}_{SU(N)}(x_1,\ldots x_\frac{N(N-1)}{2})\left(\prod_{j=1}^N\xi_j\right)\cdot\\
			&\left(\prod_{j>k}\left(\xi_j^2(1-\xi_k^2)-(1-\xi_j^2)\xi_k\right)\right)\cdot\\
			&\tilde{J}_{SU(N)}(x_{\frac{N(N-1)}{2}+1},\ldots,x_{N(N-1)}),
		\end{align*} then $\vec{0}$ does not lie in the convex hull of $\mathrm{Sp}(f)$.
	\end{conjecture}
	
	\begin{theorem}\label{thm:Mathieu_assuming_xz_conjecture_Sp(N)}
		Assume Conjecture \ref{con:xz-conjecture_Sp(N)} is true. Then the Mathieu Conjecture is true for $Sp(N)$.
	\end{theorem}
	
	\begin{conjecture}\label{con:xz-conjecture_G_2}
		Let $f:[0,1]^{6}\times \mathbb{T}^{8}\rightarrow\C$ be an admissible function in the sense of Definition \ref{def:admissible_function}. If $$\int_{\mathbb{T}^8}\int_{[0,1]^5}\int_0^{S(\xi_1)}f^P \tilde{J}_{G_2}\,d\xi_2 d\xi_1 dx_1\ldots dx_4\frac{dz_1}{z_1}\ldots\frac{dz_8}{z_8} = 0$$ for all $P\in \N$, where \begin{align}
			\begin{split}
				\tilde{J}_{G_2}(x_1,\ldots,x_4,\xi_1,\xi_2):=& \xi_1\xi_2\bigg[\xi_1^2(16(1-\xi_2^2)^3+9(1-\xi_2^2) -24(1-\xi_2^2)^2)-\\
				& (1-\xi_1^2)(3\xi_2-4\xi_2^2)^2\bigg]\big[\xi_1^2(1-\xi_2^2)-(1-\xi_1^2)\xi_2^2\big]x_1x_2x_3x_4,\end{split}
		\end{align} then $\vec{0}$ does not lie in the convex hull of $\mathrm{Sp}(f)$.
	\end{conjecture}
	
	\begin{theorem}
		Assume Conjecture \ref{con:xz-conjecture_G_2} is true. Then the Mathieu Conjecture is true for $G_2$.
	\end{theorem}

	\newpage
	\bibliography{bibfile}
	
\end{document}